\documentclass[a4paper]{elsarticle}

\usepackage[latin1]{inputenc} % LaTeX, comprends les accents !
\usepackage[T1]{fontenc}      % Police contenant les caractères français
\usepackage{geometry}         % Définir les marges
\usepackage[francais]{babel}  % Placez ici une liste de langues, la
\usepackage{amsmath,amsfonts,amssymb,epsfig}                              % dernière étant la langue principale

\newtheorem{thm}{Theorem}[section]

\newtheorem{defi}{Definition}[section]
\newtheorem{lem}{Lemma}[section]

\newtheorem{prop}{Proposition}[section]

\journal{Statistics and Probability Letters}

\begin{document}

\begin{frontmatter}



\title{A lower bound on the expected optimal value of certain random linear programs and application to shortest paths and reliability}


\author{St\'ephane Chr\'etien and Franck Corset}

\address{National Physical Laboratory, Teddington, UK and Laboratoire Jean Kuntzmann, UMR5224, Univ. Grenoble Alpes, FRANCE}

\begin{abstract}
The paper studies the expectation of the inspection time in
complex aging systems. Under reasonable assumptions, this problem
is reduced  to studying the expectation of the length of the
shortest path in the directed degradation graph of the systems
where the parameters are given by a pool of experts. The expectation 
itself being sometimes out of reach, in closed form or even through Monte
Carlo simulations in the case of large systems, we propose an easily computable lower bound. 
The proposed bound applies to a rather general class of linear programs with random nonnegative costs and is directly inspired from the upper bound of Dyer, Frieze and McDiarmid [Math.Programming {\bf 35} (1986), no.1,3--16]. .  % Résumé de l'article
\end{abstract}


\end{frontmatter}
%\title{A lower bound on the expected optimal value of certain random linear programs and application to shortest paths and reliability}           % Les paramètres du titre : titre, auteur, date
%\author{St\'ephane Chr\'etien \and Franck Corset}
% \date{}                     % La date n'est pas requise (la date du
                              % jour de compilation est utilisée en son
			      % absence

%begin{document}
%\maketitle                  % Faire un titre utilisant les données
                              % passées à \title, \author et \date

% \tableofcontents            % Table des matières
% \listoffigures              % Table des figures
% \listoftables               % Liste des tableaux
%\part{Titre}                  % Commencer une partie...
%\section{Titre}               % Commencer une section, etc.
%\subsection{Titre}            % Section plus petite
%\subsubsection{Titre}         % Encore plus petite
%\paragraph{Titre}             % Toutes petites sections (le nom \paragraph
                              % n'est pas très bien choisi)%\subparagraph{Titre}          % La dernière
%\appendix                     % Commençons les annexes
%\section{Titre}               % Annexe A
%\section{Titre}               % Annexe B

\section{Introduction and motivations}
\label{motiv}
The random shortest path problem may be a good model for describing the time to failure of 
very complex systems with various degradation schemes as for instance nuclear plants. In this section, 
we describe our motivations for studying such random shortest path problems.  

\subsection{Problem statement}
Consider a complex system whose $n$ degradation states have been
identified by experts. Let node 1 represent the state where the
system is considered as new and let node $n$ be the state of
unacceptable degradation. All maximum paths from any node of the
graph end at node $n$ as in the figure below. The system is
supposed to possibly evolve from a degradation state to any neighbor in the
corresponding connected directed acyclic graph. The transition time between any two
given states is assumed to follow a Weibull distribution whose
parameters are estimated if the number of
observations is sufficiently large. Otherwise, it is possible to make Bayesian inference in order to combine the real data with some expert opinions.

\begin{figure}
\begin{center}
\includegraphics[scale=1.2]{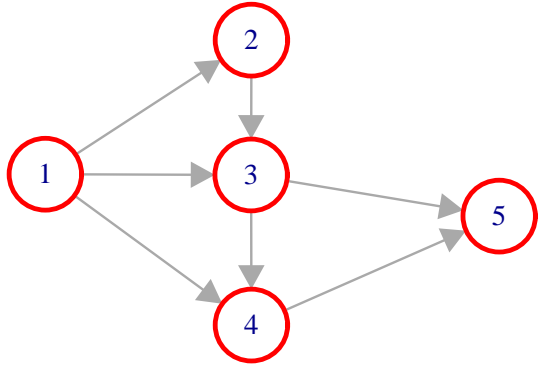}
\end{center}
%\vspace{3cm}
\label{graph}
\caption{Graph of degradation states with Weibull transitions.}
\end{figure}

Assume we start with a brand new system. Then, evolution of the
system starts in state $1$. Maintenance policies require that the
system be inspected before reaching state $n$, i.e. unacceptable
degradation. We represent this by a connected graph $\mathcal G=(V,E)$, where 
$|V|=n$ and $|E|=m$. Such examples of complex systems have 
been studied in \cite{Cor03,Cor06,Adj12}. Moreover, Chen et al. \cite{chen99} study the shortest path, in the maintenance optimization context, for some multi-state parallel-series systems. The problem posed in this paper is to provide a lower bound on acceptable inspection times.

\subsection{Inspection times and shortest paths}
In order to simplify the analysis, we assume that evolution inside
the degradation graph, a Directed Acyclic Graph (DAG), proceeds following the rule that starting
from one node $i$, the system goes to state $j$ minimizing the
transition time among neighbors of state $i$. Therefore,
acceptable inspection times will be the times lower than {\em the
shortest path} from state $1$ to state $n$ where each edge is
weighted by its transition time. In general situations, we thus
may ask for
\begin{itemize}
\item an estimator of the expected length of the shortest path from $1$ to $n$,
\item a confidence interval for the expected time path.
\end{itemize}
This task is in general impossible to achieve because of the huge number of 
observations this should require in practice. The goal of this
paper is to propose a lower bound on the expected length of the
shortest path. 

The paper is organized as follows. 
Since shortest path problems are well known to be representable as linear programs, we will address in the next section 
the more general problem of deriving a lower bound to the expectation of linear programs with random costs. In the third 
section, we specialize the study of this lower bound to an appropriate linear programming formulation of the shortest path problem. 
Moreover, we show that in the case of exponentially distributed random costs, the
Dyer-Frieze-McDiarmid upper bound is as bad as possible. The fourth section is devoted to the application to reliability 
theory as motivated by the introductory example above. In particular, the Weibull distribution is proved to satisfy the 
assumptions under which the proposed lower bound holds. 

\section{Random linear programs} 
Consider the linear program with random costs given by 
\begin{equation}
\label{randlin}
\begin{array}{rl}
z = & \min c^t x \\
& Ax=b \\
& x \geq 0 
\end{array}
\end{equation}
where $c$ is a random vector taking values on $\mathbb R_+^n$, $A$ is a matrix in $\mathbb R^{m\times n}$ and $b$ is a vector in $\mathbb R^m$.
The expectation of $c_i$ is denoted by $\mu_i$ and its variance by $\sigma_i^2$. 

In the sequel, we assume that $A$ is full rank and that the constraints of (\ref{randlin})
define a polytope which is therefore a compact set. For any subset of $\{1,\ldots,n\}$, we denote by $A_B$ the matrix whose
column set is the set of columns of $A$ indexed by $B$. We will also use the notation $x_B$ and $c_B$ for the vectors whose components 
are the components of $x$ and $c$ which are indexed by $B$. A set of indices $B$ is called a basis if its cardinality is $m$ and the matrix 
$A_B$ is full rank. A basis is said to be optimal if there exists $x^*$ in $\mathbb R^n$ such that $x^*$ is a solution to 
(\ref{randlin}) and $x^*_{B^c}=0$. 

Random linear programs have been investigated recently and many impressive 
results have been optained in the case of i.i.d. cost vectors. For 
instance, the assignment problem was investigated in \cite{wastlund2005proof}, \cite{Kro2007}, \cite{Kro2009} 
in the asymptotic regime. 

In this section, we propose a lower bound on the expected value of random linear programs in the spirit of the Dyer, Frieze and McDiarmid inequality \cite{Dyer86}. The Dyer-Frieze-McDiarmid inequality is a powerful tool for the analysis of some linear programming and combinatorial optimization problems with random costs, as detailed in the monograph of Steele \cite{Steele:SIAM97}. More precisely, 
The Dyer-Frieze-McDiarmid bounds reads as follows. 

\begin{thm}{\bf (Dyer-Frieze-McDiarmid)}
Assume that all the components of $c$ in \eqref{randlin} are independent and nonnegative and there exists $\beta \in (0,1]$
such that 
\begin{align*}
\mathbb E \left[c_i \mid c_i\ge h \right] & \ge \mathbb E \left[c_i\right]+\beta h. 
\end{align*}
Let $x$ be a feasible solution of \eqref{randlin}. Let $z^*$ denote the random optimal value of 
\eqref{randlin}. Then, assuming $\mathbb E\left[ c_1\right] x_1 \ge 
\cdots \ge \mathbb E\left[ c_n\right] x_n$. 
\begin{align*}
    \mathbb E \left[ z^*\right] & \le \beta^{-1} \sum_{i=1}^n \mathbb E\left[ c_i\right] x_i. 
\end{align*}
\end{thm}

The Weibull distributions $\mathcal W(\eta,\gamma)$ ($\eta$ and $\gamma$ are respectively the scale and shape parameters), has density function
\begin{equation}
\label{weibdf}    
f(x) = \frac{\gamma}{\eta} \left(\frac{x}{\eta}\right)^{\gamma-1} \exp\left(-\left(\frac{x}{\eta}\right)^{\gamma}\right).
\end{equation}

When the edges are Weibull distributed with shape parameters $\gamma_i$ in the interval $[1,2]$ for $i=1,\ldots,n$, we will see in Proposition \ref{MRTF}
\begin{equation}
\label{ineg}
{\rm E} [c_i \mid c_i\geq h]\leq E[c_i]+\beta h
\end{equation}
with $\beta=1$. Note that Dyer-Frieze-McDiarmid requires the reverse inequality instead, in order to hold. We will however use this property to obtain a lower bound on the expectation of the optimal value of random linear programs in Theorem \ref{main}.

As in \cite{Dyer86}, we will need the following result which is well known to users of the simplex algorithm. 
\begin{lem}
\label{lemopt}
A necessary condition for a basis to be optimal is that 
\begin{equation}
\label{opt}
c_{B^c}^t \geq c_B^t(A_B^t)^{-1}A_{B^c}. 
\end{equation}
\end{lem}

\begin{defi}
For a basis $B$, let $I_B$ be the index set 
\begin{align*}
    I_B &  = \left\{i \in \{1,\ldots,n\} \mid \left(c_B^t(A_B^t)^{-1}A_{B^c}\right)_i \ge 0 \right\}.
\end{align*}
\end{defi}

Using this result and following the same reasoning as in the proof of the Dyer, Frieze and McDiarmid inequality in \cite{Dyer86}, we obtain the following proposition.  
\begin{prop}
\label{main}
Consider the random linear program (\ref{randlin}) with random cost vector $c$ satisfying (\ref{ineg}) with $\beta \in[1,+\infty)$. Let $\mathcal B$ be a set 
of bases. Let $\mathcal I_{\mathcal B}$ be the index set 
$\mathcal I_{\mathcal B} = \cap_{B\in \mathcal B} \ I_B$.
Let $x$ be any vector satisfying the constraints of (\ref{randlin}) and such that 
\begin{align}
x_{\mathcal I_{\mathcal B}}^c & =0.    
\label{arnak}
\end{align}
Then, we have 
\begin{equation}
\label{dfmbis}
{\rm E}[z]\geq \frac1{\beta} \sum_{B\in \mathcal B}  p_B {\rm E} [c_B] x_B. 
\end{equation}
\end{prop} 
{\bf Proof}. 
Fix a basis $B$ and let $E_B$ be the event that $B$ be optimal. Take any $x$
satisfying the primal constraints. Then we have 
\begin{equation}
\begin{array}{rl}
{\rm E} [z\mid E_B] & ={\rm E} [c_B^t(A_B)^{-1}b  \mid E_B] \\
& ={\rm E} [c_B^t(A_B)^{-1}(A_Bx_B+A_{B^c}x_{B^c})  \mid E_B] \\
& ={\rm E} [c_B^tx_B+c_B^t (A_B)^{-1}A_{B^c}x_{B^c}  \mid E_B]. \\
\end{array}
\end{equation}
But using (\ref{ineg}) together with \eqref{arnak}, we have ${\rm E} [c_{B^c}^t \mid E_B,c_B]x_{B^c}\leq {\rm E} [c_{B^c}^t]x_{B^c}+ \beta c_B^t (A_B)^{-1}A_{B^c}x_{B^c}$, and thus 
\begin{equation}
\begin{array}{rl}
{\rm E} [z\mid E_B] & = {\rm E} [c_B^t\mid E_B] x_B+{\rm E} \Big[{\rm E}[ c_B^t (A_B)^{-1}A_{B^c} \mid E_B,c_B] \mid E_B\Big] x_{B^c}\\
& \geq {\rm E} [c_B^t \mid E_B]x_B +\frac1{\beta}\Big({\rm E}[c_{B^c}^t \mid E_B]-{\rm E}[c_{B^c}^t]\Big)x_{B^c} \\
& = \Big(1-\frac1{\beta}\Big){\rm E} [c_B^t \mid E_B]x_B +\frac1{\beta}{\rm E}[c^t\mid E_B]x -\frac1{\beta}{\rm E}[c_{B^c}^t]x_{B^c} .
\end{array}
\end{equation}
Since $\beta \in[1,+\infty)$ we can rule out the term $\Big(1-\frac1{\beta}\Big){\rm E} [c_B^t \mid E_B]x_B$ and using the fact that 
${\rm E}[c_{B^c}^t]x_{B^c}={\rm E}[c^t]x-{\rm E}[c_B^t]x_B$, we get 
\begin{equation}
{\rm E} [z\mid E_B]  \geq\frac1{\beta}{\rm E}[c^t\mid E_B]x -\frac1{\beta}\Big({\rm E}[c^t]x -{\rm E}[c_B^t]x_B\Big).
\end{equation}
Finally take the expectation over all possible bases to obtain 
\begin{equation}
\begin{array}{rl}
{\rm E} [z]  & \geq \frac1{\beta}\sum_B p_B {\rm E}[c^t\mid E_B]x -\frac1{\beta}\sum_B p_B {\rm E}[c^t]x+\frac1{\beta}\sum_B p_B {\rm E}[c_B^t]x_B\\
& \\
& = \frac1{\beta}\sum_B p_B {\rm E}[c_B^t]x_B.
\end{array}
\end{equation}
Restricting the summation to a special subset $\mathcal B$ of bases preserves the previous inequality and the proposition is proved. 
\hfill$\Box$

The result of this proposition is not completely satisfactory since the probabilities $p_B$ that $B$ be an optimal basis are not known. 
In certain cases, efficient approximations of these probabilities can be obtained using a more precise expression of $p_B$. Since in 
the case where the components of the cost vector $c$ are independent we easily get such an expression from the conditions for 
optimality given in Lemma \ref{lemopt}. The lower bound we thus obtain is summarized in the following theorem. 
\begin{thm}
\label{mainth}
Consider the random linear program (\ref{randlin}) with random cost vector $c$ with independent coordinates. 

a. Let $B$ be 
a basis for this program and for all $j\in B$ and $i\in B^c$, let $\alpha_{ij}=((A_B^t)^{-1}A_{B^c})_{ji}$. 
Then, we have 
\begin{equation}
p_B= {\rm E}[  \prod_{i\in B^c} P(c_i\geq \sum_{j\in B} c_j \alpha_{ji} \mid c_B)  ].    
\end{equation}

b. Let $x$ be any vector satisfying \eqref{arnak} and the constraints of (\ref{randlin}). Then
\begin{equation}
E[z]\geq \frac1{\beta}\sum_{B\in \mathcal B} {\rm E}[  \prod_{i\in B^c} P(c_i\geq \sum_{j\in B} c_j \alpha_{ji} \mid c_B)] E[c_B]^t x_B. 
\end{equation}
\end{thm}
{\bf Proof}. a. Due to independence of the components of $c$, conditionally on the value of $c_j$, $j\in B$, the events 
$c_i\geq \sum_{j\in B} c_j \alpha_{ji}$ are independent. Thus, the desired formula. 

b. Combine a. with Proposition \ref{main}. \hfill $\Box$

With these results in hand, we will now be able to turn to the more specialized case of random shortest paths in the next section.  
\section{Random shortest paths}

\subsection{Linear programming formulation}
The shortest path problem can be represented as an equivalent linear
programming problem, as is well known \cite{Hoff1963}. In \cite[pp. 75--79]{Papa98} for instance, the shortest path problem
is shown to be equivalent to
\begin{equation}
\label{linear0}
\min_{x\in \mathbb R^n} c_0^t x \text{ subject to } A_0x=b \text{ and } x\geq 0,
\end{equation}
where $c_0$ is the column vector whose components are the transition times on each edge,
$A_0$ is the incidence matrix of the oriented degradation graph and $b$ is the vector $[-1,0,\cdots,0,1]^t$,
encoding the fact that we start the path at node $1$ and end it at node $n$. Recall that the incidence
matrix is constructed as follows. Its rows are indexed by the nodes
of the graph while its columns are indexed by its edges with an extra column of all ones. In each
column indexed by edge $(i,j)$, set the $i^{\rm th}$ component to
-1, the $j^{\rm th}$ component to 1 and set all other entries to
zero. For instance, the incidence matrix for the graph of figure
\ref{graph} is given by
$$
A_0=\left[
\begin{array}{cccccccc}
-1 & -1  & -1  &  0  &  0  &  0  &  0 \\
1  &  0  &  0  & -1  &  0  &  0  &  0 \\
0  &  1  &  0  &  1  & -1  & -1  &  0 \\
0  &  0  &  1  &  0  &  1  &  0  & -1 \\
0  &  0  &  0  &  0  &  0  &  1  &  1
\end{array}
\right].
$$
Any solution vector $x^*$ to this linear program whose components are binary, i.e. $\in \{0,1\}$ encodes a path
whose edges correspond to the nonzero components of $x^*$. The important property is that the matrix $A_0$ is totally unimodular (TUM) 
which means that every square submatrix has determinant equal to -1, 0 or 1. This linear programming formulation of the problem has however a drawback: the incidence matrix is not full rank and the size of its kernel is the number of connected components of the graph \cite{Bollobas:GTM98}. On the other hand, for our results to apply recall that we need the matrix $A$
in (\ref{randlin}) to be full rank. In order to remedy this problem, we introduce the extended incidence matrix $A$, given by
$$ A=[A_0 \mid e],$$
where $e$ is the vector whose components are all equal to one. For instance, the extended incidence matrix for the graph
of figure \ref{graph} is given by
$$
A=\left[
\begin{array}{cccccccc}
-1 & -1  & -1  &  0  &  0  &  0  &  0 & 1 \\
1  &  0  &  0  & -1  &  0  &  0  &  0 & 1 \\
0  &  1  &  0  &  1  & -1  & -1  &  0 & 1 \\
0  &  0  &  1  &  0  &  1  &  0  & -1 & 1 \\
0  &  0  &  0  &  0  &  0  &  1  &  1 & 1
\end{array}
\right].
$$
In addition, let $c$ denote the extended cost vector $[c_0^t,0]$. Then, we get the following proposition.
\begin{prop}
\label{extend}
The shortest path problem is equivalent to the linear program
\begin{equation}
\label{linear}
\min_{x\in \mathbb R^n} c^t x \text{ subject to } Ax=b \text{ and } x\geq 0.
\end{equation}
\end{prop}
{\bf Proof}. Let $x$ be an optimal solution of the given linear program. Then, $x$ satisfies the Karush-Kuhn-Tucker equations which are of the form:

\begin{equation}
\label{kkt}
\left \lbrace \begin{array}{l}
 c+
\left[
\begin{array}{cc}
I & A^t
\end{array}
\right]
\left[
\begin{array}{c}
u \\
v
\end{array}
\right] =0 ,\\
Ax=b \\
x\geq 0 \\
u_i x_i=0 \text{ for } i=1,\ldots,n \\
u_i \leq 0 \\
v_j(a_j^t x-b_j)=0 \text{ for } j=1,\ldots,m
\end{array}
\right.
\end{equation}
where $a_j^t$ is the $j^{th}$ row of $A$ and the vectors $u\in \mathbb R^n$ and $v\in \mathbb R^m$ are the Lagrange multipliers.
More precisely, the multipliers that compose the vector $u$ deal with the nonegativity constraints and the components of $v$ deal
with the others. The third equation is imposed in order to select the "active" constraints at optimality.
In particular, it implies that if $u_{n+1}\neq 0$ we must have $x_{n+1}=0$. On the other hand, if $u_{n+1}=0$,
$x_{n+1}$ may be positive. In what follows, we show that $x_{n+1}$ is always null which will readily imply that
this linear program also solves the shortest path problem.

Equations (\ref{kkt}) determine a polyhedron in $\mathbb R^{n+1}\times \mathbb R^{n}\times \mathbb R^{m}$.
Now assume that $[x^t,u^t,v^t]^t$ is a corner point of this polyhedron and that $u_{n+1}=0$. Since the matrix
$A$ is now full rank, the $x$-part of the corner vector satisfies
$$ A x=b, $$
and
$$ x_i=0, \hspace{.3cm} i\in I_x $$
for some index set $I_x$ of cardinality $n+1-m$. Now write these last nullity constraints $Cx=0$ for some appropriate matrix $C$.
Then by Cramer's rules, we obtain that the last coordinate $x_{n+1}$ is proportional to 
$$
{\rm det}
\Big(\left[
\begin{array}{c}
C \\
\left[
\begin{array}{cc}
A_0 \mid b
\end{array}
\right]
\end{array}
\right]
\Big).
$$
On the other hand, we know that the sum of the rows of $A$ is equal to zero and the same holds for the sum of
the components of $b$. Therefore, the determinant just above is null. Therefore $x_{n+1}=0$ as announced.
From this, it is easy to deduce that the vector of the first $n$ components of an optimal solution
to the present linear program also solves (\ref{linear0}).
\hfill$\Box$

Using this proposition, we deduce that Theorem \ref{mainth} applies to the random shortest path problem. We now consider the Dyer-Frieze-McDiarmid bound which gives an upper bound on the 
expected optimal value.  

\subsection{The DFM upper bound}
The upper bound for the expected optimal cost of random linear programs Dyer, Frieze and Mc Diarmid 
\cite{Dyer86} is a very nice result and a major contribution to the study of random optimization problems;
see also \cite{Steele:SIAM97}. It can be stated as follows. 
\begin{thm}{\rm (\cite{Dyer86})}
Assume that the random costs $c_i$ are independent and satisfy
$$
E[c_i\mid c_i\geq h]\geq E[c_i]+\alpha h
$$
for some $\alpha \in (0,1]$. \footnote{This property is obviously satistied with equality and $\alpha=1$ in the case of
exponential transitions.} Then for any matrix $A\in \mathbb R^{n\times m}$ and any vector
$b\in \mathbb R^n$, the optimal value $z^*$ of the general linear program (\ref{linear}) satisfies
\begin{equation}
\label{dfm}
E[z^*]\leq \max_{S \: :\: \# S=n}\sum_{i\in S}E[c_i]x_i
\end{equation}
for any feasible solution $x$, i.e. any $x$ satisfying $Ax=b$.
\end{thm} 

It is interesting to understand to what extent the Dyer-Frieze-McDiarmid (DFM) bound is useful for the shortest 
path problem. Surprisingly, the answer is that the DFM bound is as bad as possible in this case, despite is remarkable 
efficiency on other standard combinatorial problems as shown in \cite[Chapter 4]{Steele:SIAM97}. To understand 
why this happens, consider the deterministic problem where the random costs are replaced by their expected values. 
\begin{equation}
\begin{array}{rl}
\label{minexp}
\zeta =  & \min  E[c]^T x \\
&  Ax=b\\
&  x \geq 0
\end{array}
\end{equation}
Now, notice that whatever the distribution of the cost vector $c$ may be, the following upper bound is immediate to obtain:
\begin{equation}
{\rm E}[z]\leq \zeta.
\end{equation}
The following proposition shows that the DFM bound is no better than this trivial upper bound. 
\begin{prop}
Consider problem (\ref{randlin}) where the random costs are assumed
to be independent and exponentially distributed. Then the DFM bound is equal to the optimum 
value $\zeta$ of the associated deterministic program (\ref{minexp}). 
\end{prop}
{\bf Proof}. Take $x$ equal to any binary
vector minimizing (\ref{minexp}). It is clear that the number of
ones in this vector is less than the number of nodes in the graph.
Then, the maximum value over all sets $S$ of cardinality $n$ in
the right hand term in (\ref{dfm}) is obtained when $S$ is taken
to be the set of indices $i$ for which $x_i=1$. Thus $\sum_{i\in S} E[c_i] x_i$ is exactly the cost of $x$, i.e. $\zeta$.
\hfill$\Box$

Contrary to intuition, replacing the
random costs by their expected values is far from being a safe idea for the problem of providing efficient lower bounds 
to the mean inspection time.

\section{Application to reliability}
In this section, we address the problem of finding lower bounds to the inspection time of complex systems 
in reliability. As explained in Section \ref{motiv} our main interest in studying random shortest paths problems 
relies in its possible application to the analysis of the time to failure for very complex systems. We will 
assume in this section that the transition times between two degradation states follows a Weibull distribution.    
In order to apply our previous results, we will need to study the Weibull distribution a little further.  
\subsection{Some properties of the Weibull distribution}
Let $X$ be a random variable with Weibull distribution $Weib(\eta,\gamma)$, i.e. with probability density function given by 
\begin{equation}
f_X(t)=\frac{\gamma}{\eta}\Big(\frac{t}{\eta} \Big)^\gamma
e^{-(\frac{t}{\eta})^\gamma}.
\end{equation}
Then, the mean residual time to failure (MRTF) is given by
\begin{equation}
G_X(h)= E[X\mid X\geq h]=\eta
e^{(\frac{h}{\eta})^\gamma} \Gamma(1+\frac{1}{\gamma},(\frac{h}{\eta})^\gamma),
\end{equation}
where $\Gamma(a,h)$ is the incomplete gamma function defined by
\begin{equation}
\Gamma(a,h)=\int_{h}^{+\infty} t^{a-1} e^{-t} dt.
\end{equation}
\begin{lem}
\label{propWeib}
Let $X$ be a Weibull $Weib(\eta,\gamma)$ distributed random variable. Then, 

a. the first two derivatives of the MRTF for a Weibull distributed
variable $X$ are given by
\begin{equation}
\label{Gprime}
G_X^{\prime}(h)=\gamma \Big(\frac{h}{\eta}\Big)^{\gamma}\Big( \frac{\eta}{h} e^{(\frac{h}{\eta})^\gamma}\Gamma\big( 1+\frac{1}{\gamma},(\frac{h}{\eta})^\beta\big) -1\Big) 
\end{equation}
and
\begin{eqnarray}
\begin{array}{rll}
G_X^{\prime\prime}(h) & = & -\frac{\gamma^2}{\eta^\gamma}h^{(\gamma-1)}\Big( 1+\gamma(\frac{h}{\eta})^\gamma\Big) +\\
 &  & \frac{\gamma}{\eta^{(\gamma-1)}}h^{(\gamma-2)}e^{(\frac{h}{\gamma})^\gamma}\Gamma(1+\frac{1}{\gamma},(\frac{h}{\eta})^\gamma)\Big( \gamma(\frac{h}{\eta})^\gamma+\gamma-1\Big).
\end{array}
\end{eqnarray}
Moreover, 

b. when $\gamma\geq 1$ we have
\begin{equation}
\lim_{h\to 0} G_X^\prime(h)=0, \hspace{.3cm} \lim_{h\to +\infty} G_X^\prime(h)=1 \text{ and }  G_X^{\prime\prime}(h)>0.
\end{equation}
\end{lem}
{\bf Proof}. a. We omit the proof of the formula for the first and second derivative of $G_X$. 

b. Now, since $G_X^\prime$ is clearly continuous on $\mathbb R^+$ and $G_X^\prime(0)=0$, we obtain the first 
assertion in b. Using the transformation $u=t-(h/\eta)^\gamma$, the first of the two terms between parenthesis in (\ref{Gprime}) can we written 
\begin{equation}
\eta e^{\displaystyle(\frac{h}{\eta})^\gamma}\Gamma\big( 1+\frac{1}{\gamma},(\frac{h}{\eta})^\gamma\big)=h\int_0^{+\infty}e^{-u}\big( 1+(\frac{\eta}{h})^\gamma u\big )^{\frac{1}{\gamma}}\,du.
\end{equation}
For all $u\geq 0$ we have the following Taylor expansion 
\begin{equation}
( 1+(\frac{\eta}{h})^\gamma u\big)^{\frac{1}{\gamma}}=1+\frac{1}{\gamma}\frac{\eta^\gamma}{h^\gamma}u+o(\frac{1}{h^{\gamma}}). 
\end{equation}
Multiplying by $e^{-u}$ and integrating, we obtain
\begin{equation}
\begin{array}{rl}
\int_0^{+\infty}e^{-u}\big( 1+(\frac{\eta}{h})^\gamma u\big )^{\frac{1}{\gamma}}\,du  & =\displaystyle{\int_0^{+\infty}e^{-u}\,du+\frac{1}{\gamma}\frac{\eta^\gamma}{h^\gamma}\int_0^{+\infty}e^{-u}u\,du} + o(\frac{1}{h^{\gamma}})\\
& \\
& = \displaystyle{1+\frac{1}{\gamma}\frac{\eta^\gamma}{h^\gamma} +o(\frac{1}{h^{\gamma}})}.
\end{array}
\end{equation}
Combining with (\ref{Gprime}), we deduce that $\lim_{h\rightarrow +\infty} G_X^{\prime}(h)=1$. Finally, it is readily seen that $G_X^{\prime\prime}$ is positive as soon as $\gamma >1$.  
\hfill$\Box$

\subsection{The lower bound}
In this subsection, we work out an easily computable lower bound derived from Theorem \ref{mainth}. We first have 
the following crucial result saying that the most commonly encountered Weibull distributions in 
reliability theory satisfy the main assumption of Proposition \ref{main} and Theorem \ref{mainth}.
\begin{prop}
\label{MRTF}
Assume that $X$ has distribution $Weib(\eta,\gamma)$ with $\gamma \in (1,2)$. Then, for all $h\geq 0$, we have
$$
E[X\mid X\geq h] \leq E[X]+h.
$$
\end{prop}
{\bf Proof}.
This is a direct consequence of Proposition \ref{propWeib}. 
\hfill$\Box$   

In the next theorem, we derive an explicit lower bound from Theorem \ref{mainth} in the case of Weibull distributions. 
\begin{thm}
\label{mainth2}
Consider the random linear program (\ref{randlin}) with random cost vector $c$ with independent components and 
assume that each component $c_i$, $1\leq i\leq n$ follows a Weibull distribution $Weib(\eta_i,\gamma_i)$ 

a. Let $B$ be a basis for this program and for all $j\in B$ and $i\in B^c$, let $\alpha_{ij}=((A_B^t)^{-1}A_{B^c})_{ji}$. 
Then, we have 
\begin{equation}
p_B\geq 1-\sum_{i\in B^c}\frac{4 \sqrt{e} \Big(\sum_{j\in B}\alpha_{ji}^2\eta_j^2 \Gamma \Big(1+\frac{2}{\gamma_j}\Big)\Big)^{\frac12} +|\sum_{j\in B}\alpha_{ji}\eta_j \Gamma \Big(1+\frac{1}{\gamma_j}\Big) |}{\eta_i^{\gamma_i}}.    
\end{equation}

b. Let $x$ be any vector satisfying \eqref{arnak} and the constraints of (\ref{randlin}). Then
\begin{equation}
E[z]\geq \sum_{B\in \mathcal B} \Big(1-\sum_{i\in B^c}\frac{4 \sqrt{e} \Big(\sum_{j\in B}\alpha_{ji}^2\eta_j^2 \Gamma \Big(1+\frac{2}{\gamma_j}\Big)\Big)^{\frac12} +|\sum_{j\in B}\alpha_{ji}\eta_j \Gamma \Big(1+\frac{1}{\gamma_j}\Big) |}{\eta_i^{\gamma_i}}\Big) E[c_B]^t x_B. 
\end{equation}
\end{thm}

{\bf Proof}. a. Theorem \ref{mainth}.a. gives the following formula for $p_B$: 
\begin{equation}
p_B= {\rm E}\Big[ \prod_{i \in B^c} \exp -\Big(\max \Big\{0,\frac{\sum_{j\in B}\alpha_{ji}c_j}{\eta_i} \Big\} \Big)^{\gamma_i}\Big], 
\end{equation}
where the expectation is taken with respect to the variables $c_j$, $j\in B$. Now since $\exp (-x)\geq 1-x$, we obtain that 
\begin{equation}
\begin{array}{rl}
p_B & \geq  {\rm E} \Big[1-\sum_{i \in B^c}\max \Big\{0,\Big(\frac{\sum_{j\in B}\alpha_{ji}c_j}{\eta_i} \Big)^{\gamma_i} \Big\}\Big], \\ 
& = 1-\sum_{i\in B^c} {\rm E}\Big[ \max \Big\{0,\frac{\sum_{j\in B}\alpha_{ji}c_j}{\eta_i} \Big\}^{\gamma_i}\Big].
\end{array}
\end{equation}
Thus, 
\begin{equation}
p_B  \geq  1-\sum_{i \in B^c} \frac{{\rm E} \Big[\max\Big\{0,\sum_{j\in B}\alpha_{ji}c_j \Big\}^{\gamma_i}\Big]}{\eta_i^{\gamma_i}}. 
\end{equation}
In order to simplify the subsequent computations, we will use the crude majorization: $$\max\Big\{0,\sum_{j\in B}\alpha_{ji}c_j \Big\}\leq |\sum_{j\in B}\alpha_{ji}c_j|$$ 
which gives 
\begin{equation}
\label{tata}
p_B  \geq  1-\sum_{i \in B^c} \frac{{\rm E} \Big[|\sum_{j\in B}\alpha_{ji}c_j |^{\gamma_i}\Big]}{\eta_i^{\gamma_i}}. 
\end{equation}
Now, our next goal is to use the Kintchine inequalities in order to bound the last expression by a quantity expressed in terms of the 
$l_2$ norm which will be easier to control. For this purpose, one might want to center the random variables involved in (\ref{tata}) and use the triangle inequality to obtain 
\begin{equation}
\label{titi}
\begin{array}{rl}
{\rm E} \Big[|\sum_{j\in B}\alpha_{ji}c_j |^{\gamma_i}\Big]^{\frac1{\gamma_i}} & = {\rm E} \Big[|\sum_{j\in B}\alpha_{ji}(c_j-{\rm E}[c_j]+{\rm E}[c_j]) |^{\gamma_i}\Big]^{\frac1{\gamma_i}} \\
& \leq {\rm E} \Big[|\sum_{j\in B}\alpha_{ji}(c_j-{\rm E}[c_j]) |^{\gamma_i}\Big]^{\frac1{\gamma_i}}+{\rm E} \Big[|\sum_{j\in B}\alpha_{ji}{\rm E}[c_j] |^{\gamma_i}\Big]^{\frac1{\gamma_i}} \\
& = {\rm E} \Big[|\sum_{j\in B}\alpha_{ji}(c_j-{\rm E}[c_j]) |^{\gamma_i}\Big]^{\frac1{\gamma_i}}+|\sum_{j\in B}\alpha_{ji}{\rm E}[c_j] |
\end{array}
\end{equation}
Using Jensen's inequality, a standard trick gives 
\begin{equation}
\begin{array}{rl}
{\rm E} \Big[|\sum_{j\in B}\alpha_{ji}(c_j-{\rm E}[c_j]) |^{\gamma_i}\Big]^{\frac1{\gamma_i}} & = {\rm E} \Big[|\sum_{j\in B}\alpha_{ji}(c_j-{\rm E}[c_j^\prime]) |^{\gamma_i}\Big]^{\frac1{\gamma_i}}\\
& \leq {\rm E} \Big[|\sum_{j\in B}\alpha_{ji}(c_j-c_j^\prime) |^{\gamma_i}\Big]^{\frac1{\gamma_i}}
\end{array}
\end{equation}
where $c_j^\prime$, $j\in B$ are i.i.d. variables independent of $c_j$, $j\in B$ and such that $c_j$ has same distribution as $c_j^\prime$, $j\in B$. Let 
$\epsilon_j$, $j\in B$ be standard Rademacher $\pm 1$ random variables. Since $\sum_{j\in B}\alpha_{ji}(c_j-c_j^\prime)$ has the same distribution as 
$\sum_{j\in B}\alpha_{ji}\epsilon_j (c_j-c_j^\prime)$, we have 
\begin{equation}
\label{tutu}
\begin{array}{rl}
{\rm E} \Big[|\sum_{j\in B}\alpha_{ji}(c_j-{\rm E}[c_j]) |^{\gamma_i}\Big]^{\frac1{\gamma_i}} 
& = {\rm E} \Big[|\sum_{j\in B}\alpha_{ji}\epsilon_j (c_j-c_j^\prime) |^{\gamma_i}\Big]^{\frac1{\gamma_i}}\\
& \leq 2 {\rm E} \Big[|\sum_{j\in B}\alpha_{ji}\epsilon_j c_j |^{\gamma_i}\Big]^{\frac1{\gamma_i}}
\end{array}
\end{equation}
where we used once again the triangle inequality. Notice that  
\begin{equation}
\label{toto}
\begin{array}{rl}
{\rm E} \Big[|\sum_{j\in B}\alpha_{ji}\epsilon_j c_j |^{\gamma_i}\Big]^{\frac1{\gamma_i}} & = {\rm E} \Big[{\rm E} \Big[|\sum_{j\in B}\alpha_{ji}\epsilon_j c_j |^{\gamma_i}\mid c_j,\: j\in B \Big]\Big]^{\frac1{\gamma_i}} \\
& = {\rm E} \Big[ \Big({\rm E} \Big[|\sum_{j\in B}\alpha_{ji}\epsilon_j c_j |^{\gamma_i}\mid c_j,\: j\in B \Big]^{\frac1{\gamma_i}}\Big)^{\gamma_i}\Big]^{\frac1{\gamma_i}}. 
\end{array}
\end{equation}
On the other hand, Khintchine's inequality gives 
\begin{equation}
\label{khin}
\begin{array}{rl}
{\rm E} \Big[\Big({\rm E} \Big[|\sum_{j\in B}\alpha_{ji}\epsilon_j c_j |^{\gamma_i}\mid c_j,\: j\in B \Big]^{\frac1{\gamma_i}}\Big)^{\gamma_i}\Big]^{\frac1{\gamma_i}}
& \leq C_{\gamma_i}{\rm E} \Big[\Big((\sum_{j\in B}\alpha_{ji}^2c_j^2 )^{\frac12}\Big)^{\gamma_i}\Big]^{\frac1{\gamma_i}} 
\end{array}
\end{equation}
where $C_{\gamma_i}$ is equal to $\sqrt{2e\gamma_i}$ in the present context. Thus, 
\begin{equation}
\begin{array}{rl}
{\rm E} \Big[\Big({\rm E} \Big[|\sum_{j\in B}\alpha_{ji}\epsilon_j c_j |^{\gamma_i}\mid c_j,\: j\in B \Big]^{\frac1{\gamma_i}}\Big)^{\gamma_i}\Big]^{\frac1{\gamma_i}}
& \leq \sqrt{2e\gamma_i}\Big({\rm E} \Big[(\sum_{j\in B}\alpha_{ji}^2c_j^2 )^{\frac{\gamma_i}2}\Big]^{\frac2{\gamma_i}}\Big)^{\frac12}. 
\end{array}
\end{equation}
Moreover, since $p\mapsto E[|X|^p]^{\frac1p}$ is an increasing function and $\gamma_i$ is assumed to belong to $[1,2]$, we obtain the simpler bound
\begin{equation}
\begin{array}{rl}
{\rm E} \Big[\Big({\rm E} \Big[|\sum_{j\in B}\alpha_{ji}\epsilon_j c_j |^{\gamma_i}\mid c_j,\: j\in B \Big]^{\frac1{\gamma_i}}\Big)^{\gamma_i}\Big]^{\frac1{\gamma_i}}
& \leq 2\sqrt{e}{\rm E} \Big[\sum_{j\in B}\alpha_{ji}^2c_j^2 \Big]^{\frac12},\\
& \\
& = 2\sqrt{e}\Big(\sum_{j\in B}\alpha_{ji}^2{\rm E}\big[c_j^2\big]\Big)^{\frac12} 
\end{array}
\end{equation}
Moreover, since ${\rm E} [c_j^2]= \eta_j^2 \Gamma \Big(1+\frac{2}{\gamma_j}\Big)$,
we have 
\begin{equation}
\begin{array}{rl}
{\rm E} \Big[\Big({\rm E} \Big[|\sum_{j\in B}\alpha_{ji}\epsilon_j c_j |^{\gamma_i}\mid c_j,\: j\in B \Big]^{\frac1{\gamma_i}}\Big)^{\gamma_i}\Big]^{\frac1{\gamma_i}}
& \leq 2\sqrt{e}{\rm E} \Big[\sum_{j\in B}\alpha_{ji}^2c_j^2 \Big]^{\frac12},\\
& \\
& = 2\sqrt{e}\Big(\sum_{j\in B}\alpha_{ji}^2\eta_j^2 \Gamma \Big(1+\frac{2}{\gamma_j}\Big)\Big)^{\frac12}. 
\end{array}
\end{equation}
Combining this result with (\ref{tata}), (\ref{titi}), (\ref{tutu}), (\ref{toto}), and replacing 
${\rm E} [c_j]= \eta_j \Gamma \Big(1+\frac{1}{\gamma_j}\Big)$ in (\ref{titi}),
we finally obtain the desired result. 

b. This follows from part a. and Proposition \ref{main}. 
\hfill$\Box$

\section{Conclusion and perspectives}

In this paper, we derived a lower bound on the probability that a given path is optimal for the shortest path problem with 
independent arc weights with Weibull distributions. For this purpose, we used the linear programming formulation of the problem
and extended the work of Dyer, Frieze and Mc Diarmid \cite{Dyer86}. The results presented here are of a theoretical 
nature. Further refinements and applications to real data will be proposed in a subsequent paper.

\section*{References}

\bibliographystyle{plain}
\bibliography{DFM}

\end{document}